\documentclass[times]{iapress}
\usepackage{moreverb}

\usepackage[colorlinks,bookmarksopen,bookmarksnumbered,citecolor=red,urlcolor=red]{hyperref} 

\def\volumeyear{2019}
\def\volumenumber{7}
\def\volumemonth{March}
\usepackage{amsmath,amssymb}
\usepackage{graphicx}

\newcommand{\lp}[1]{\left( \begin{array}{#1} }
\newcommand{\rp}{\end{array} \right)}

\newcommand{\be}{\begin{equation}}
\newcommand{\ee}{\end{equation}}
\newtheorem{thm}{Theorem}[section]
\newtheorem{lem}{Lemma}[section]
\newtheorem{nas}{Corrolary}[section]
\theoremstyle{Remark}

\theoremstyle{definition}
\newtheorem{ozn}{Definition}[section]

\begin{document}

\runningheads{Interpolation of Multidimensional Process with Missing Observations}{O. Masyutka, M. Moklyachuk  and M. Sidei}

\title{Interpolation   Problem  for Multidimensional Stationary Processes with Missing Observations}

\author{Oleksandr Masyutka\affil{1}, Mikhail Moklyachuk \affil{2}$^,$\corrauth,
Maria Sidei\affil{2}}

\address{
\affilnum{1} Department of Mathematics and Theoretical Radiophysics,
Taras Shevchenko National University of Kyiv, Kyiv 01601, Ukraine
\affilnum{2}Department of Probability Theory, Statistics and Actuarial
Mathematics, Taras Shevchenko National University of Kyiv, Kyiv 01601, Ukraine
}
\corraddr{Mikhail Moklyachuk (Email: Moklyachuk@gmail.com). Department of Probability Theory, Statistics and Actuarial
Mathematics, Taras Shevchenko National University of Kyiv, Volodymyrska 64 Str., Kyiv 01601, Ukraine.}

 \begin{abstract}
The  problem  of the mean-square optimal linear estimation of linear
functionals  which depend on the unknown values of a multidimensional continuous time  stationary stochastic process is considered.
Estimates are based on observations of the process with an additive stationary stochastic noise process at points which do not belong to some finite intervals of a real line.
The problem is investigated in the case of spectral certainty, where the spectral densities of the processes are exactly known.
Formulas for calculating the mean-square errors and spectral characteristics of the optimal linear estimates of functionals are proposed under the condition of spectral certainty.
The minimax (robust) method of estimation is applied in the case spectral uncertainty, where spectral densities of the processes
are not known exactly while some sets of admissible spectral densities of the processes are given.
Formulas that determine the least favorable spectral
densities and the minimax spectral characteristics of the optimal estimates of functionals are proposed for some special sets of admissible spectral densities.
\end{abstract}

\keywords{Stationary process, optimal estimate, mean square error, minimax-robust estimate, least favorable spectral density, minimax spectral characteristic}

\maketitle
\noindent{\bf AMS 2010 subject classifications.} Primary: 60G10, 60G25, 60G35, Secondary: 62M20, 93E10, 93E11

\vspace{10pt}

\noindent {\bf DOI:} 10.19139/soic.v7i1.430

\section{Introduction}
 The problem of estimation of the unknown values of stochastic processes is of constant interest in the theory and applications of stochastic processes.
 The description of the problem and effective methods of solution of the interpolation, extrapolation and filtering problems for stationary stochastic sequences and processes were
 were developed by Kolmogorov \cite{Kolmogorov}, Wiener \cite{Wiener} and Yaglom \cite{Yaglom1,Yaglom2}.
Further development of the methods of solution and analysis of  the problems can be found in the works by  Rozanov \cite{Rozanov} and Hannan \cite{Hannan},
 Box et. al \cite{Box:Jenkins},  Brockwell and  Davis \cite{Brockwell_Davis}.
 The basic assumption of most of the methods of solution of the problem of estimation of the unobserved values of stochastic processes is that the spectral densities of the considered stochastic processes are exactly known.
However, in practice these methods can not be applied since the complete information on the spectral densities is impossible in most cases.
In order to solve the problem, the parametric or nonparametric estimates of the unknown spectral densities are found.
Then, one of the traditional estimation methods is applied, provided that the selected spectral densities are true.
This procedure can result in significant increase in the value of the error of estimate as Vastola and Poor \cite{Vastola} have demonstrated with the help of some examples.
To avoid this effect one can search estimates which are optimal for all densities from a certain class of admissible spectral densities.
These estimates are called minimax since they minimize the maximum value of the error.
The paper by Grenander \cite{Grenander} should be marked as the first one where this approach to extrapolation problem for stationary processes was proposed.
Several models of spectral uncertainty and the minimax-robust methods of data processing can be found in the survey paper by Kassam and Poor \cite{Kassam}.
In their papers, Franke \cite{Franke1,Franke},  Franke and Poor \cite{Franke_Poor} presented results of the investigation of the minimax extrapolation  and filtering problems for stationary sequences with the help of convex optimization methods. This approach makes it possible to find relations that determine the least favorable spectral densities for different classes of admissible spectral densities.
The papers  by Moklyachuk \cite{Moklyachuk:2008,Moklyachuk:2015} is dedicated to the investigation of the problems of the  optimal linear estimation of  functionals which depend on the unknown values of stationary sequences and processes.
The problem of estimation of functionals which depend on the unknown values of multivariate stationary stochastic
processes is the aim of the papers by Moklyachuk and Masyutka \cite{Moklyachuk:Mas2006}--\cite{Moklyachuk:Mas2012}.
In the book  by Moklyachuk and Golichenko \cite{Golichenko} results of investigation of the interpolation, extrapolation and filtering problems for periodically correlated
stochastic processes are described.
 Estimation problems for functionals which depend on the unknown values of the stochastic processes with stationary increments were investigated by Luz and Moklyachuk  \cite{luz7}--\cite{Luz2017}.
The problem of interpolation, extrapolation and filtering for stationary processes with missing observations were investigated by Moklyachuk and Sidei \cite{Sidei1} -- \cite{Sidei7}.
The prediction and interpolation problems for stationary stochastic sequences with missing observations were investigated in papers by Bondon \cite{Bondon1, Bondon2},
Cheng, Miamee and Pourahmadi \cite{Cheng1},
Cheng and Pourahmadi \cite{Cheng2},
Kasahara,  Pourahmadi and Inoue \cite{Kasahara},
Pourahmadi, Inoue and Kasahara \cite{Pourahmadi},
Pelagatti \cite{Pelagatti}. The problem of interpolation of  stationary sequences was considered in the paper of Salehi \cite{Salehi}.

In this paper we deal with the problem of the mean-square optimal linear estimation of the functional $$A_s\vec{\xi}=\sum\limits_{l=1}^{s}\int\limits_ {-M_{l}-N_{l}}^{-M_{l}}\vec{a}(t)^\top\vec{\xi}(t)dt,$$ which depends on the unknown values of a multidimensional continuous time  stationary stochastic process
 $\vec{\xi}(t)$
  from observations of the process with an additive stationary stochastic noise process $\vec{\xi}(t)+\vec{\eta}(t)$ at points  $t\in\mathbb{R}  \backslash S$,
$S=\bigcup\limits_{l=1}^{s}[ -M_{l}-N_{l}, -M_{l} ]$. The case of spectral certainty as well as the case  of spectral uncertainty are considered.
In the case of spectral certainty, where the spectral densities of the processes are exactly known, formulas for calculating the spectral characteristic and the mean square error of the optimal linear estimate of the functional are derived.
In the case of spectral uncertainty, where the spectral densities are not exactly known while a set of admissible spectral densities is given, the minimax method is applied.
This method allows us to find estimates that minimize the maximum values of the mean-square errors of the estimates for all spectral density matrices from a given class of admissible spectral density matrices.
Formulas that determine the least favorable spectral densities and the minimax-robust spectral characteristics of the optimal estimates of the functional are proposed for some specific classes of admissible spectral densities.

\section{Hilbert space projection approach to interpolation problem for stationary processes with missing observations}

Let $\vec{\xi}(t)=\left \{ \xi_ {k} (t) \right \}_{k = 1} ^ {T}, t\in \mathbb{R},$ and $\vec{\eta}(t)=\left \{ \eta_ {k} (t) \right \}_{k = 1} ^ {T}, t\in  \mathbb{R},$ be uncorrelated mean square continuous multidimensional stationary stochastic processes with zero first moments, $E\vec{\xi}(t)=\vec{0}$, $E\vec{\eta}(t)=\vec{0}$, and
 correlation functions
which admit the spectral decomposition (see Gikhman and Skorokhod \cite{Gihman})
$$
 R_{\xi}(t)=\,\frac{1}{2\pi}\int\limits_{-\infty}^{\infty}e^{it\lambda}F(\lambda)d\lambda, \quad R_{\eta}(t)=\,\frac{1}{2\pi}\int\limits_{-\infty}^{\infty}e^{it\lambda}G(\lambda)d\lambda,
 $$
where $F(\lambda)=\left\{f_{kl}(\lambda)\right\}_{k,l=1}^T$ and $G(\lambda)=\left\{g_{kl}(\lambda)\right\}_{k,l=1}^T$ are spectral density matrices of the processes $\vec{\xi}(t)$ and $\vec{\eta}(t)$ respectively which satisfy the minimality condition
\begin{equation}\label{minimal}
\int _{-\infty }^{\infty } (b(\lambda ))^{\top} (F(\lambda )+G(\lambda
))^{-1} \overline{b(\lambda )}d\lambda < \infty,
\end{equation}
where $b(\lambda)=\sum\limits_{l=1}^{s}\int\limits_ {-M_{l}-N_{l}}^{-M_{l}}\vec{\alpha}(t)e^{it\lambda}dt$ is a nontrivial function of the exponential type.
Under this condition the error-free interpolation of the process $\vec{\xi}(t)+\vec{\eta}(t)$ is impossible  (see Rozanov \cite{Rozanov}).

The stationary processes $\vec{\xi}(t)$ and $\vec{\eta}(t)$ admit the spectral decompositions (see Gikhman and Skorokhod \cite{Gihman}; Karhunen \cite{Karhunen})
\begin{equation} \label{ksi}
\vec{\xi}(t)=\int\limits_{-\infty}^{\infty}e^{it\lambda}Z_{\xi}(d\lambda), \hspace{1cm}
\vec{\eta}(t)=\int\limits_{-\infty}^{\infty}e^{it\lambda}Z_{\eta}(d\lambda),
\end{equation}
where $Z_{\xi}(d\lambda)$ and $Z_{\eta}(d\lambda)$ are vector valued orthogonal stochastic measures such that
\[
EZ_{\xi}(\Delta_1)(Z_{\xi}(\Delta_2))^*=\,\frac{1}{2\pi}\int_{\Delta_1\cap\Delta_2}F(\lambda)d\lambda,
\]
\[
EZ_{\eta}(\Delta_1)(Z_{\eta}(\Delta_2))^*=\,\frac{1}{2\pi}\int_{\Delta_1\cap\Delta_2}G(\lambda)d\lambda.
\]

Consider the problem of the mean-square optimal linear estimation of the functional
$$A_s\vec{\xi}=\sum\limits_{l=1}^{s}\int\limits_ {-M_{l}-N_{l}}^{-M_{l}}\vec{a}(t)^\top\vec{\xi}(t)dt,$$
 which depends on the unknown values of the process $\vec{\xi}(t)$. Estimates are based on observations of the process $\vec{\xi}(t)+\vec{\eta}(t)$ at points  $t\in\mathbb{R}  \backslash S$, where $S=\bigcup\limits_{l=1}^{s}[ -M_{l}-N_{l}, -M_{l} ]$.

We will assume that the function $\vec{a}(t)$ satisfies the following conditions
\begin{equation}
\sum\limits_{k=1}^{T}\sum\limits_{l=1}^{s}\int\limits_ {-M_{l}-N_{l}}^{-M_{l}}|{a_k}(t)|dt<\infty,\quad
\sum\limits_{l=1}^{s}\int\limits_ {-M_{l}-N_{l}}^{-M_{l}}t\|\vec{a}(t)\|^2dt<\infty.
\end{equation}
This condition ensures that the functional  $A_s\vec{\xi}$ has a finite second moment.

Making use of the spectral decomposition (\ref{ksi}) of the stationary process $\vec{\xi}(t)$, we can represent the functional $A_s\vec{\xi}$ in the following form
\begin{equation*}
A_s\vec{\xi}=\int\limits_{-\infty}^{\infty}(A_s(\lambda))^\top Z_{\xi}(d\lambda), \quad
  A_s(\lambda)=\sum\limits_{l=1}^{s}\int\limits_ {-M_{l}-N_{l}}^{-M_{l}}\vec{a}(t)e^{it\lambda }dt.
 \end{equation*}

Denote by $\hat{A}_s\vec{\xi}$ the optimal linear estimate of the functional $A_s\vec{\xi}$ from  observations of the process $\vec{\xi}(t)+\vec{\eta}(t)$
and denote by $\Delta(F,G)=E\left|A_s\vec{\xi}-\hat{A}_s\vec{\xi}\right|^2$ the mean-square error of the estimate $\hat{A}_s\vec{\xi}$.
Since the spectral densities of the stationary processes $\vec{\xi}(t)$ and $\vec{\eta}(t)$ are suppose to be known, we can use the method of orthogonal projections in Hilbert spaces proposed by Kolmogorov \cite{Kolmogorov} to find the estimate $\hat{A}_s\vec{\xi}$.

Consider   values $\xi_k(t),k=1,\dots,T,t\in \mathbb{R}$ and $\eta_k(t),k=1,\dots,T,t\in \mathbb{R}$ as   elements of the Hilbert space $H=L_2(\Omega,\mathcal{F},P)$ generated by random variables $\xi$ with zero mathematical expectations, $E\xi=0$,  finite variations, $E|\xi|^2<\infty$, and inner product $(\xi,\eta)=E\xi\overline{\eta}$.
Denote by $H^s(\xi+\eta)$ the closed linear subspace generated by elements $\{\xi_k(t)+\eta_k(t): t\in \mathbb{R}  \backslash S, k=\overline{1,T}\}$ in the Hilbert space $H=L_2(\Omega,\mathcal{F},P)$.

Denote by $L_{2} (F+G)$ the Hilbert space of vector-valued functions $ \vec{a}( \lambda )= \left \{a_{k} ( \lambda ) \right \}_{k=1}^{T} $ such that
 \[ \int_{- \infty}^{ \infty} \vec{a}( \lambda )^{ \top} (F( \lambda )+G( \lambda )) \overline{ \vec{a}( \lambda )}d \lambda  = \int_{- \infty}^{ \infty} \, \sum_{k,l=1}^{T}a_{k} ( \lambda ) \overline{a_{l} ( \lambda )} \, (f_{kl} ( \lambda )+g_{kl} ( \lambda )) \, d \lambda  < \infty, \]
and denote by $L_2^s(F+G)$  the subspace of the space  $L_2(F+G)$ generated by the functions $\{e^{it\lambda } \delta _{k},\; \delta _{k} =\left\{\delta _{kl} \right\}_{l=1}^{T} , k=\overline{1,T}, t\in \mathbb{R} \backslash S \}.$

The mean-square optimal linear estimate $\hat{A}_s\vec{\xi}$ of the functional  $A_s\vec{\xi}$ can be represented in the form
 \begin {equation*}
\hat{A}_s\vec{\xi}=\int\limits_{-\infty}^{\infty}(h(\lambda))^\top(Z_{\xi}(d\lambda)+ Z_{\eta}(d\lambda)),
 \end{equation*}
where $h(\lambda) \in L_2^s(F+G)$ is the spectral characteristic of the estimate.

The mean-square error $\Delta(h;F,G)$ of the estimate $\hat{A}_s\vec{\xi}$ is given by the formula
\begin{equation*}
 \Delta(h;F,G)=E\left|A_s\xi-\hat{A}_s\xi\right|^2=\]
\[=\frac{1}{2\pi}\int\limits_{-\infty}^{\infty}( A_{s} (\lambda )-h(\lambda ))^{\top} F(\lambda )\overline{(A_{s} (\lambda )-h(\lambda ))}d\lambda
 +\frac{1}{2\pi}\int\limits_{-\infty}^{\infty}(h(\lambda ))^{\top}
G(\lambda )\overline{h(\lambda )}d\lambda.
\end{equation*}

According to the Hilbert space projection method  the optimal estimation of the functional $A_s\vec{\xi}$ is a projection of the element $A_s\vec{\xi}$ of the space $H$  on  the subspace $H^s(\xi+\eta)$. It can be found from the following conditions:
\begin{equation*} \begin{split}
1)& \hat{A}_s\vec{\xi} \in H^s(\xi+\eta), \\
2)& A_s\vec{\xi}-\hat{A}_s\vec{\xi} \bot  H^s(\xi+\eta).
\end{split} \end{equation*}

It follows from the second condition that the spectral characteristic  $h(\lambda)$ of the optimal linear estimate $\hat{A}_s\vec{\xi}$ for any $t\in \mathbb{R} \backslash S $ satisfies the equation
\begin{equation*}
\frac{1}{2\pi}\int\limits_{-\infty}^{\infty} \left(A_s(\lambda)- h(\lambda)\right)^\top F(\lambda)e^{-it\lambda}d\lambda-\frac{1}{2\pi}\int\limits_{-\infty}^{\infty}  (h(\lambda))^\top G(\lambda)e^{-it\lambda}d\lambda=0.
\end{equation*}

\noindent The last relation is equivalent to equations
\begin{equation}\label{cond2} \begin{split}
\frac{1}{2\pi}\int\limits_{-\infty}^{\infty} \left[(A_s(\lambda))^\top F(\lambda)- (h(\lambda))^\top (F(\lambda)+G(\lambda)))\right]e^{-it\lambda}d\lambda=0, \quad  t\in \mathbb{R} \backslash S.
\end{split} \end{equation}

\noindent Define the function
$$(C_s(\lambda))^\top=(A_s(\lambda))^\top F(\lambda)- (h(\lambda))^\top (F(\lambda)  +G(\lambda))$$
and its Fourier transform
$$\vec{c}(t)=\frac{1}{2\pi}\int\limits_{-\infty}^{\infty} C_s(\lambda)e^{-it\lambda}d\lambda,\quad t \in \mathbb{R}.$$

\noindent It follows from  the  relation (\ref{cond2}) that the function $\vec{c}(t)$ is nonzero only on the set $S$. Hence, the function $C_s(\lambda)$ is of the form
$$C_s(\lambda)=\sum\limits_{l=1}^{s}\int\limits_ {-M_{l}-N_{l}}^{-M_{l}} \vec{c}(t)e^{it\lambda}dt, $$
 and the spectral characteristic of the estimate $\hat{A}_s\vec{\xi}$ is of the form
\begin{equation} \label{sphar} \begin{split}
(h(\lambda))^\top=(A_s(\lambda))^\top F(\lambda)(F(\lambda)+G(\lambda))^{-1}-(C_s(\lambda))^\top (F(\lambda)+G(\lambda))^{-1}.
\end{split} \end{equation}

It follows from the first condition, $\hat{A}_s\vec{\xi} \in H^s(\xi+\eta)$, which determine the optimal linear estimate of the functional   $A_s\vec{\xi}$, that
the following relation holds true
\begin{equation}\label{13} \begin{split}
\int\limits_{-\infty}^{\infty}\left((A_s(\lambda))^\top F(\lambda)(F(\lambda)+G(\lambda))^{-1}-(C_s(\lambda))^\top (F(\lambda)+G(\lambda))^{-1}\right)e^{-it\lambda}d\lambda=0,\, t\in  S.
\end{split} \end{equation}

The last relation can be represented in terms of linear operators in the space $L_2(S)$. Let us define  operators
\begin{equation*}\begin{split}
(\bold{B}_s \bold{a})(t)&=\frac{1}{2\pi}\sum\limits_{l=1}^{s}\int\limits_ {-M_{l}-N_{l}}^{-M_{l}}\int\limits_{-\infty}^{\infty}(\vec{a}(u))^\top(F(\lambda)+G(\lambda))^{-1}e^{i\lambda(u-t)}d\lambda du, \\
(\bold{R}_s \bold{a})(t)&=\frac{1}{2\pi}\sum\limits_{l=1}^{s}\int\limits_ {-M_{l}-N_{l}}^{-M_{l}}\int\limits_{-\infty}^{\infty}(\vec{a}(u))^\top F(\lambda)(F(\lambda)+G(\lambda))^{-1}e^{i\lambda(u-t)}d\lambda du,\\
(\bold{Q}_s \bold{a})(t)&=\frac{1}{2\pi}\sum\limits_{l=1}^{s}\int\limits_ {-M_{l}-N_{l}}^{-M_{l}}\int\limits_{-\infty}^{\infty}(\vec{a}(u))^\top F(\lambda)(F(\lambda)+G(\lambda))^{-1}G(\lambda)e^{i\lambda(u-t)}d\lambda du,\\
&\vec{a}(t) \in L_2(S), \quad t \in S.
\end{split}\end{equation*}

Making use of the introduced operators,
 relation (\ref{13}) can be represented in the form
\begin{equation}\label{rivn2}
(\bold{R}_s\bold{a})(t)=(\bold{B}_s\bold{c})(t), \quad t \in S.
\end{equation}

Suppose that  the operator $\bold{B}_s$ is invertible (see paper by Salehi \cite{Salehi} for more details). Then the function $\vec{c}(t)$ can be calculated by the formula
$$\vec{c}(t)=(\bold{B}_s^{-1}\bold{R}_s\bold{a})(t), \quad t \in S.$$

Consequently, the spectral characteristic $h(\lambda)$ of the estimate $\hat{A}_s\vec{\xi}$ is calculated by the formula
 \begin{equation}\label{4} \begin{split}
&(h(\lambda))^\top=(A_s(\lambda))^\top F(\lambda)(F(\lambda)+G(\lambda))^{-1}-
(C_s(\lambda))^\top(F(\lambda)+G(\lambda))^{-1}, \\
&C_s(\lambda)=\sum\limits_{l=1}^{s}\int\limits_ {-M_{l}-N_{l}}^{-M_{l}} (\bold{B}_s^{-1}\bold{R}_s\bold{a})(t)e^{it\lambda}dt, \quad   t \in S.
\end{split} \end{equation}

The mean-square error of the estimate  $\hat{A}_s\vec{\xi}$   can be calculated by the formula
\[\Delta(h;F,G)=\frac{1}{2\pi}\int\limits_{-\infty}^{\infty}((A_{s} (\lambda ))^{\top} G(\lambda )+(C_{s} (\lambda ))^{\top} ) (F(\lambda )+G(\lambda ))^{-1}  F(\lambda )\times\]
\[\times(F(\lambda )+G(\lambda ))^{-1}((A_{s} (\lambda ))^{\top} G(\lambda )+(C_{s} (\lambda ))^{\top} )^{*}d\lambda+\]
\[+\frac{1}{2\pi}\int\limits_{-\infty}^{\infty}((A_{s} (\lambda ))^{\top} F(\lambda )-(C_{s} (\lambda ))^{\top} )(F(\lambda )+G(\lambda ))^{-1} G(\lambda )\times\]
\[\times (F(\lambda )+G(\lambda ))^{-1}((A_{s} (\lambda ))^{\top} G(\lambda )+(C_{s} (\lambda ))^{\top} )^{*}d\lambda=\]
\begin{equation} \label{55}
=\langle(\bold{R}_s\bold{a})(t),(\bold{B}_s^{-1}\bold{R}_s\bold{a})(t)\rangle+\langle(\bold{Q}_s\bold{a})(t),\vec{a}(t)\rangle,
\end{equation}
where $\left\langle a(t),b(t)\right\rangle $ is the inner product in the space $L_2(S)$.

Let us summarize the obtained results and present them in the form of a theorem.

\begin{thm} \label{t1}
Let $\{\vec{\xi}(t), t\in \mathbb{R}\}$ and $\{\vec{\eta}(t), t\in  \mathbb{R}\}$ be uncorrelated multidimensional stationary stochastic processes with the spectral density matrices $F(\lambda)$ and $G(\lambda)$ which satisfy the minimality condition (\ref{minimal}).
The spectral characteristic   $h(\lambda)$ and the mean-square error  $\Delta(F,G)$ of the optimal linear estimate of the functional $A_s\vec{\xi}$ which depends on the unknown values of the process  $\vec{\xi}(t)$ based on observations of the process  $\vec{\xi}(t)+\vec{\eta}(t),$ $t\in \mathbb{R}\backslash S$ can be calculated by formulas (\ref{4}), (\ref{55}).
\end{thm}

Consider the case where the stationary process $\vec{\xi}(t)$ is observed without noise. In this case the spectral characteristic of the estimate $\hat{A}_s\vec{\xi}$ is of the form
 \begin{equation} \label{sphar1} \begin{split}
&(h(\lambda))^\top=(A_s(\lambda))^\top-(C_s(\lambda))^\top (F(\lambda))^{-1},\\
&C_s(\lambda)=\sum\limits_{l=1}^{s}\int\limits_ {-M_{l}-N_{l}}^{-M_{l}} \vec{c}(t)e^{it\lambda}dt,
\end{split} \end{equation}
the relation (\ref{rivn2}) is of the form
\begin{equation} \label{rivn}
   \vec{a}(t)=(\bold{B}_s\bold{c})(t), \quad t \in S.
 \end{equation}

Suppose that the operator $\bold{B}_s$ is invertible. Then the unknown function $\vec{c}(t)$ can be calculated by the formula
$$\vec{c}(t)=(\bold{B}_s^{-1}\bold{a})(t), \quad t \in S,$$
and the spectral characteristic of the estimate $\hat{A}_s\vec{\xi}$ is of the form
\begin {equation} \label{spchar2} \begin{split}
&(h(\lambda))^\top=(A_s(\lambda))^\top-(C_s(\lambda))^\top (F(\lambda))^{-1}, \\
&C_s(\lambda)=\sum\limits_{l=1}^{s}\int\limits_ {-M_{l}-N_{l}}^{-M_{l}} (\bold{B}_s^{-1}\bold{a})(t)e^{it\lambda}dt, \quad   t \in S.
\end{split} \end{equation}

The mean-square error of the estimate can be calculated by the formula
\begin{equation} \label {6} \begin{split}
 \Delta(h;F)=\left\langle(\bold{B}_s^{-1})\bold{a}(t), \vec{a}(t)\right\rangle.
\end{split} \end{equation}

The following theorem holds true.

\begin{thm} \label{t2}
Let $\{\vec{\xi}(t), t\in \mathbb{R}\}$  be a multidimensional stationary stochastic process with the spectral density matrix $F(\lambda),$ which satisfies the minimality condition
\begin{equation}\label{minimal1}
\int\limits_{-\infty}^{\infty}(b(\lambda ))^{\top} (F(\lambda ))^{-1} \overline{b(\lambda )}d\lambda<\infty
\end{equation}
for a nonzero vector-valued function of the exponential type $b(\lambda)=\sum\limits_{l=1}^{s}\int\limits_ {-M_{l}-N_{l}}^{-M_{l}}\vec{\alpha}(t)e^{it\lambda}dt$.
The spectral characteristic $h(\lambda)$ and the mean-square error  $\Delta(F)$ of the optimal linear estimate of the functional $A_s\vec{\xi}$ which depends on the unknown values of the process  $\vec{\xi}(t)$ based on observations of the process $\vec{\xi}(t)$ at points $t\in\mathbb{R} \backslash S $,  $S=\bigcup\limits_{l=1}^{s}[-M_{l}-N_l, -M_{l} ],$ can be calculated by formulas (\ref{spchar2}), (\ref{6}).
\end{thm}


\section{Minimax approach to interpolation problem for stationary processes with missing observations}

The derived formulas can be applied only in the case of spectral certainty, where the spectral density matrices of the considered processes are exactly known.
However, in practice,  we do not have exact values of  the spectral density matrices while, instead, we have  a  class of admissible spectral density matrices, where the spectral density matrices of the processes belong to. In this case, the minimax method can be applied to estimate the value of the functional.
This method allows us to find estimates that minimize the maximum values of the mean-square errors of the estimates for all spectral density matrices from a given class of admissible spectral density matrices.
For description of the minimax method we introduce the following definitions (see
Moklyachuk \cite{Moklyachuk:2008, Moklyachuk:2015}, and Moklyachuk and Masytka \cite{Moklyachuk:Mas2012}).

 \begin{ozn}
 For a given class of spectral densities $D=D_F \times D_G$ the spectral densities  $F^0(\lambda) \in D_F$, $G^0(\lambda) \in D_G$ are called the least favorable in the class $D$ for the optimal linear interpolation of the functional $A_s\vec{\xi}$  if the following relation holds true $$\Delta\left(F^0,G^0\right)=\Delta\left(h\left(F^0,G^0\right);F^0,G^0\right)=\max\limits_{(F,G)\in D_F\times D_G}\Delta\left(h\left(F,G\right);F,G\right).$$
\end{ozn}

\begin{ozn}
For a given class of spectral densities $D=D_F \times D_G$ the spectral characteristic $h^0(\lambda)$ of the optimal linear interpolation  of the functional $A_s\vec{\xi}$ is called minimax-robust if there are satisfied conditions
$$h^0(\lambda)\in H_D= \bigcap\limits_{(F,G)\in D_F\times D_G} L_2^s(F+G),$$
$$\min\limits_{h\in H_D}\max\limits_{(F,G)\in D}\Delta\left(h;F,G\right)=\max\limits_{(F,G)\in D}\Delta\left(h^0;F,G\right).$$
\end{ozn}

From the introduced definitions and formulas derived in previous section we can obtain the following statements.

\begin{lem}
The spectral densities $F_0(\lambda)\in D_F,$ $G_0(\lambda) \in D_G$ satisfying the minimality condition (\ref{minimal}) are the least favorable in the class  $D=D_F\times D_G$ for the optimal linear interpolation of the functional $A_s\vec{\xi},$ if the Fourier coefficients of the functions $$(F^0(\lambda)+G^0(\lambda))^{-1},\quad F^0(\lambda)(F^0(\lambda)+G^0(\lambda))^{-1}, \quad F^0(\lambda)(F^0(\lambda)+G^0(\lambda))^{-1}G^0(\lambda)$$ determine  operators $\bold{B}_s^0, \bold{R}_s^0, \bold{Q}_s^0$, which determine a solution of the constrained optimization problem
\begin{equation} \label{extrem} \begin{split}
&\max\limits_{(F,G)\in D_F\times D_G}\langle(\bold{R}_s\bold{a})(t),(\bold{B}_s^{-1}\bold{R}_s\bold{a})(t)\rangle +\langle(\bold{Q}_s\bold{a})(t),\vec{a}(t)\rangle= \\
&=\langle(\bold{R}_s^0\bold{a})(t),((\bold{B}_s^0)^{-1}\bold{R}_s^0\bold{a})(t)\rangle+\langle(\bold{Q}_s^0\bold{a})(t),\vec{a}(t)\rangle.
\end{split}\end{equation}
The minimax spectral characteristic  $h^0=h(F^0,G^0)$ is determined by the formula (\ref{4}) if $h(F^0,G^0) \in H_D.$
\end{lem}

\begin{nas}
Let the spectral density $F^0(\lambda)\in D_F$  satisfy the minimality condition (\ref{minimal1}).  The spectral density $F^0(\lambda)\in D_F$ is the least favorable in the class $ D_F$ for the optimal linear interpolation of the functional $A_s\vec{\xi}$  from  observations of the process $\vec{\xi}(t)$ at points $t \in \mathbb{R}\backslash S$ if the Fourier coefficients of the function $(F^0(\lambda))^{-1}$ determine the operator  $\bold{B}_s^0$ which determines a solution of the constrained optimization problem
 \begin{equation} \label{extrem2}
\max\limits_{F\in D_F}\langle(\bold{B}_s^{-1}\bold{a})(t),\vec{a}(t)\rangle=\langle((\bold{B}_s^0)^{-1}\bold{a})(t),\vec{a}(t)\rangle.
\end{equation}
The minimax spectral characteristic $h^0=h(F^0)$ is determined by the formula  (\ref{spchar2}) if $h(F^0) \in H_{D_F}.$
 \end{nas}

For more detailed analysis of properties of the least favorable spectral densities and the minimax-robust spectral characteristics we observe that the least favorable spectral densities
 $F^0(\lambda)$, $G^0(\lambda)$ and the minimax spectral characteristic $h^0=h(F^0,G^0)$ form a saddle point of the function $\Delta \left(h;F,G\right)$ on the set  $H_D\times D.$
 The saddle point inequalities
$$\Delta\left(h^0;F,G\right)  \leq\Delta\left(h^0;F^0,G^0\right)\leq \Delta\left(h;F^0,G^0\right), \quad \forall h \in H_D, \forall F \in D_F, \forall G \in D_G,$$
hold true if  $h^0=h(F^0,G^0)$ and $h(F^0,G^0)\in H_D,$ where $(F^0,G^0)$ is a solution to the constrained optimization problem
\begin{equation} \label{7}
\sup\limits_{(F,G)\in D_F\times D_G}\Delta\left(h(F^0,G^0);F,G\right)=\Delta\left(h(F^0,G^0);F^0,G^0\right).
\end{equation}
The linear functional $\Delta \left(h\left(F^{0} , G^{0} \right); F, G\right)$
is calculated by the formula
\[\Delta \left(h\left(F^{0} , G^{0} \right); F, G\right)=\]
\[=\frac{1}{2\pi } \int _{-\infty }^{\infty }  ((A_{s} (\lambda ))^{\top} G^{0} (\lambda )+(C_{s}^{0} (\lambda ))^{\top} ) (F^{0} (\lambda )+G^{0} (\lambda ))^{-1} F(\lambda )\times\]
\[\times (F^{0} (\lambda )+G^{0} (\lambda ))^{-1}((A_{s} (\lambda ))^{\top} G^{0} (\lambda )+(C_{s}^{0} (\lambda ))^{\top} )^{*}d\lambda+\]
\[+\frac{1}{2\pi } \int _{-\infty }^{\infty }  ((A_{s} (\lambda ))^{\top} F^{0} (\lambda )-(C_{s}^{0} (\lambda ))^{\top} )(F^{0} (\lambda )+G^{0} (\lambda ))^{-1}G(\lambda )\times\]
\[\times (F^{0} (\lambda )+G^{0} (\lambda ))^{-1} ((A_{s} (\lambda ))^{\top} G^{0} (\lambda )-(C_{s}^{0} (\lambda ))^{\top} )^{*}d\lambda,\]
\begin{equation*}
C_s^0(\lambda)=\sum\limits_{l=1}^{s}\int\limits_ {-M_{l}-N_{l}}^{-M_{l}} ((\bold{B}_s^0)^{-1}\bold{R}_s^0\bold{a})(t)e^{it\lambda}dt, \quad   t \in S.
\end{equation*}

The constrained optimization problem (\ref{7}) is equivalent to the unconstrained optimization problem (see Pshenichnyj \cite{Pshenychnyj})
\begin{equation} \label{8}
\Delta_D(F,G)=-\Delta(h(F^0,G^0);F,G)+\delta((F,G)\left|D_F\times D_G\right.)\rightarrow \inf,
\end{equation}
where $\delta((F,G)|D_F\times D_G)$ is the indicator function of the set $D=D_F\times D_G$.
A solution of the problem (\ref{8}) is characterized by the condition $0 \in \partial\Delta_D(F^0,G^0),$ where $\partial\Delta_D(F^0,G^0)$ is the subdifferential of the convex functional $\Delta_D(F,G)$ at point $(F^0,G^0)$,
which is the necessary and sufficient condition under which the pair $(F^0,G^0)$ belongs to the set of minimums of the convex functional $\Delta(h(F^0,G^0);F,G)$.
This condition makes it possible to find the least favourable spectral densities in some special classes of spectral densities $D$ (see books by Ioffe and Tihomirov \cite{Ioffe}, Pshenichnyj \cite{Pshenychnyj}, Rockafellar \cite{Rockafellar}).

Note, that the form of the functional $\Delta(h(F^0,G^0);F,G)$ is convenient for application of the Lagrange method of indefinite multipliers for finding solution to the problem (\ref{8}).
Making use of the method of Lagrange multipliers and the form of
subdifferentials of the indicator functions we describe relations that determine the least favourable spectral densities in some special classes of spectral densities (see books by Moklyachuk \cite{Moklyachuk:2008b,Moklyachuk:2008}, Moklyachuk and Masyutka \cite{Moklyachuk:Mas2012} for additional details).

Taking into consideration the introduced definitions and the derived relations we can verify that the following lemma holds true.

\begin{lem}
Let $(F^0,G^0)$ be a solution of the optimization problem (\ref{8}). The spectral densities  $F^0(\lambda)$, $G^0(\lambda)$ are the least favorable in the class $D=D_F\times D_G$ and the spectral characteristic  $h^0=h(F^0,G^0)$ is the minimax of the optimal linear estimate of the functional  $A_s\vec{\xi}$ if  $h(F^0,G^0) \in H_D$.
\end{lem}

\section{Least favorable spectral densities in the class $D=D_0 \times D_{\varepsilon}$}

Consider the problem of  mean-square optimal interpolation of the functional $A_s\vec{\xi}$  in the case where the spectral density matrices of the processes are not exactly known while the
 admissible spectral density matrices are from the class $D=D_0 \times D_{\varepsilon}$, where
$$ D_{0}^{1} =\bigg\{F(\lambda )\left|\frac{1}{2\pi }
\int _{-\infty }^{\infty }{\rm{Tr}}\, F(\lambda )d\lambda =p\right.\bigg\},$$
$$D_{\varepsilon }^{1}  =\bigg\{G(\lambda )\bigg|{\mathrm{Tr}}\,
G(\lambda )=(1-\varepsilon ) {\mathrm{Tr}}\,  G_{1} (\lambda
)+\varepsilon {\mathrm{Tr}}\,  W(\lambda ), \frac{1}{2\pi } \int _{-\infty }^{\infty }{\mathrm{Tr}}\,
G(\lambda )d\lambda =q \bigg\};$$
$$D_{0}^{2} =\bigg\{F(\lambda )\left|\frac{1}{2\pi }
\int _{-\infty }^{\infty }f_{kk} (\lambda )d\lambda =p_{k}, k=\overline{1,T}\right.\bigg\},$$
$$D_{\varepsilon }^{2}  =\bigg\{G(\lambda )\bigg|g_{kk} (\lambda)
=(1-\varepsilon )g_{kk}^{1} (\lambda )+\varepsilon w_{kk}(\lambda), \frac{1}{2\pi } \int _{-\infty }^{\infty }g_{kk} (\lambda)d\lambda  =q_{k} , k=\overline{1,T}\bigg\};$$
$$D_{0}^{3} =\bigg\{F(\lambda )\left|\frac{1}{2\pi } \int _{-\infty }^{\infty }\left\langle B_{1} ,F(\lambda )\right\rangle d\lambda  =p\right.\bigg\},$$
$$D_{\varepsilon }^{3} =\bigg\{G(\lambda )\bigg|\left\langle
B_{2},G(\lambda )\right\rangle =(1-\varepsilon )\left\langle B_{2}
,G_{1} (\lambda )\right\rangle +  \varepsilon \left\langle B_{2}
,W(\lambda )\right\rangle, \frac{1}{2\pi } \int _{-\infty }^{\infty
}\left\langle B_{2} ,G(\lambda )\right\rangle d\lambda =q\bigg\};$$
$$D_{0}^{4} =\bigg\{F(\lambda )\left|\frac{1}{2\pi } \int
_{-\infty }^{\infty }F(\lambda )d\lambda  =P\right.\bigg\},$$
$$D_{\varepsilon }^{4}=\left\{G(\lambda )\bigg|G(\lambda
)=(1-\varepsilon )G_{1} (\lambda )+\varepsilon W(\lambda ),
\frac{1}{2\pi } \int _{-\infty }^{\infty }G(\lambda )d\lambda
=Q\right\}.$$
Here  $G_1(\lambda)$ is a known and fixed spectral density matrix while $W(\lambda)$ is an unknown spectral density matrix,
$p, q, p_k, q_k, k=\overline{1,T}$ are given numbers, $P, Q, B_1, B_2$ are given positive definite Hermitian matrices.
The class $D_{\varepsilon}$ describes ``$\varepsilon$-contamination''\, model of stochastic processes.

From the condition $0\in \partial \Delta _{D} (F^{0} ,G^{0} )$ we find the following equations which determine the least favourable spectral densities for these given sets of admissible spectral densities.

For the first pair $D_{0}^{1}\times D_{\varepsilon}^{1}$ we have equations
\begin{equation} \label{eq_4_1}
((A_{s} (\lambda ))^{\top} G^{0} (\lambda )+(C_{s}^{0} (\lambda
))^{\top} )^{*}((A_{s} (\lambda ))^{\top} G^{0} (\lambda )+(C_{s}^{0}
(\lambda ))^{\top} )=\alpha ^{2} (F^{0} (\lambda )+G^{0} (\lambda
))^{2} ,
\end{equation}
\begin{equation} \label{eq_4_2}
((A_{s} (\lambda ))^{\top} F^{0} (\lambda )-(C_{s}^{0} (\lambda
))^{\top} )^{*}((A_{s} (\lambda ))^{\top} F^{0} (\lambda )-(C_{s}^{0}
(\lambda ))^{\top} )=(\beta ^{2} +\gamma(\lambda ))(F^{0} (\lambda
)+G^{0} (\lambda ))^{2},
\end{equation}
where $\alpha^{2}, \beta^{2}$ are Lagrange multipliers, $\gamma(\lambda )\le 0$ and $\gamma(\lambda )=0$ if ${\mathrm{Tr}}\,
F^{0} (\lambda )>(1-\varepsilon ) {\mathrm{Tr}}\,  G_{1} (\lambda ).$

For the second pair $D_{0}^{2}\times D_{\varepsilon}^{2}$ we have equations
\[((A_{s} (\lambda ))^{\top} G^{0} (\lambda )+(C_{s}^{0} (\lambda
))^{\top} )^{*}((A_{s} (\lambda ))^{\top} G^{0} (\lambda )+(C_{s}^{0}
(\lambda ))^{\top} )=\]
\begin{equation}  \label{eq_4_3}
=(F^{0} (\lambda )+G^{0} (\lambda
))\left\{\alpha _{k}^{2} \delta _{kl} \right\}_{k,l=1}^{T} (F^{0}
(\lambda )+G^{0} (\lambda )),
\end{equation}
\[((A_{s} (\lambda ))^{\top} F^{0} (\lambda )-(C_{s}^{0} (\lambda
))^{\top} )^{*}((A_{s} (\lambda ))^{\top} F^{0} (\lambda )-(C_{s}^{0}
(\lambda ))^{\top})=\]
\begin{equation} \label{eq_4_4}
=(F^{0} (\lambda )+G^{0} (\lambda )) \left\{(\beta _{k}^{2} +\gamma_{k} (\lambda ))\delta _{kl} \right\}_{k,l=1}^{T} (F^{0} (\lambda
)+G^{0} (\lambda )),
\end{equation}
where $\alpha _{k}^{2}, \beta_{k}^{2}$ are Lagrange multipliers, $\delta _{kl}$ are Kronecker symbols, $\gamma _{k}(\lambda )\le 0$ and $\gamma _{k} (\lambda )=0$ if $g_{kk}^{0}
(\lambda )>(1-\varepsilon )g_{kk}^{1} (\lambda ).$

For the third pair $D_{0}^{3}\times D_{\varepsilon}^{3}$ we have equations
\[((A_{s} (\lambda ))^{\top} G^{0} (\lambda )+(C_{s}^{0} (\lambda
))^{\top} )^{*}((A_{s} (\lambda ))^{\top} G^{0} (\lambda )+(C_{s}^{0}
(\lambda ))^{\top} )=\]
\begin{equation}  \label{eq_4_5}
=\alpha ^{2} (F^{0} (\lambda )+G^{0} (\lambda
))B_{1}^{\top} (F^{0} (\lambda )+G^{0} (\lambda )),
\end{equation}
\[((A_{s} (\lambda ))^{\top} F^{0} (\lambda )-(C_{s}^{0} (\lambda
))^{\top} )^{*}((A_{s} (\lambda ))^{\top} F^{0} (\lambda )-(C_{s}^{0}
(\lambda ))^{\top})=\]
\begin{equation} \label{eq_4_6}
=(\beta ^{2} +\gamma'(\lambda ))(F^{0} (\lambda )+G^{0} (\lambda ))B_{2}^{\top} (F^{0}
(\lambda )+G^{0} (\lambda )),
\end{equation}
where $\alpha^{2}, \beta^{2}$ are Lagrange multipliers, $\gamma' ( \lambda )\le 0$ and $\gamma' ( \lambda )=0$ if $\langle
B_{2} ,G^{0} ( \lambda ) \rangle
>(1- \varepsilon ) \langle B_{2} ,G_{1} ( \lambda ) \rangle.$

For the fourth pair $D_{0}^{4}\times D_{\varepsilon}^{4}$ we have equations
\[((A_{s} (\lambda ))^{\top} G^{0} (\lambda )+(C_{s}^{0} (\lambda
))^{\top} )^{*}((A_{s} (\lambda ))^{\top} G^{0} (\lambda )+(C_{s}^{0}
(\lambda ))^{\top} )=\]
\begin{equation} \label{eq_4_7}
=(F^{0} (\lambda )+G^{0} (\lambda
))\vec{\alpha }\cdot \vec{\alpha }^{*}(F^{0} (\lambda )+G^{0}
(\lambda )),
\end{equation}
\[((A_{s} (\lambda ))^{\top} F^{0} (\lambda )-(C_{s}^{0} (\lambda
))^{\top} )^{*}((A_{s} (\lambda ))^{\top} F^{0} (\lambda )-(C_{s}^{0}
(\lambda ))^{\top})=\]
\begin{equation} \label{eq_4_8}
=(F^{0} (\lambda )+G^{0}(\lambda ))(\vec{\beta }\cdot \vec{\beta }^{*}+\Gamma(\lambda
))(F^{0} (\lambda )+G^{0} (\lambda )),
\end{equation}
where $\vec{\alpha}, \vec{\beta}$ are Lagrange multipliers, $\Gamma(\lambda )\le 0$ and $\Gamma _{3} (\lambda )=0$ if $G^{0}
(\lambda )>(1-\varepsilon )G_{1} (\lambda ).$

Thus, the following statement holds true.

\begin{thm}
The least favorable spectral densities  $F^0(\lambda)$, $G^0(\lambda)$  in the classes $D_0 \times D_{\varepsilon}$ for the optimal linear interpolation of the functional $A_s\vec{\xi}$ are determined by relations
(\ref{eq_4_1}), (\ref{eq_4_2}) for the first pair $D_{0}^{1}\times D_{\varepsilon}^{1}$ of sets of admissible spectral densities;
(\ref{eq_4_3}), (\ref{eq_4_4}) for the second pair $D_{0}^{2}\times D_{\varepsilon}^{2}$ of sets of admissible spectral densities;
(\ref{eq_4_5}), (\ref{eq_4_6}) for the third pair $D_{0}^{3}\times D_{\varepsilon}^{3}$ of sets of admissible spectral densities;
(\ref{eq_4_7}), (\ref{eq_4_8}) for the fourth pair $D_{0}^{4}\times D_{\varepsilon}^{4}$ of sets of admissible spectral densities;
 the minimality condition (\ref{minimal}); the constrained optimization problem (\ref{extrem}) and restrictions  on densities from the corresponding classes $D_0 \times D_{\varepsilon}$.
 The minimax-robust spectral characteristic of the optimal estimate of the functional $A_s\vec{\xi}$ is determined by the formula (\ref{4}).
\end{thm}

\begin{nas}
 The least favorable spectral densities $F^{0}(\lambda)$ in the classes $D_0^{k}$, $k=1,2,3,4$, for the optimal linear interpolation of the functional  $A_s\vec{\xi}$, which depends on the unknown values of the process  $\vec{\xi}(t)$ based on observations of the process $\vec{\xi}(t)$ at points $t\in\mathbb{R} \backslash S$, are determined by the following  equations, respectively,
\begin{equation}
((C_{s}^{0}(\lambda) )^{\top} )^{*}\cdot(C_{s}^{0}(\lambda) )^{\top}=\alpha^{2}(F^{0} (\lambda ))^{2},
\end{equation}
\begin{equation}
((C_{s}^{0}(\lambda) )^{\top} )^{*}\cdot(C_{s}^{0}(\lambda) )^{\top}=F^{0} (\lambda )\left\{\alpha _{k}^{2} \delta _{kl} \right\}_{k,l=1}^{T} F^{0} (\lambda ),
\end{equation}
\begin{equation}
((C_{s}^{0}(\lambda) )^{\top} )^{*}\cdot(C_{s}^{0}(\lambda) )^{\top}=\alpha^{2}F^{0} (\lambda )(B_2)^\top F^{0} (\lambda ),
\end{equation}
\begin{equation}
((C_{s}^{0}(\lambda) )^{\top} )^{*}\cdot(C_{s}^{0}(\lambda) )^{\top}=F^{0} (\lambda )\vec{\alpha}\cdot \vec{\alpha}^{*}F^{0} (\lambda ),
\end{equation}
 the minimality condition (\ref{minimal1}); the constrained optimization problem (\ref{extrem2}) and restrictions  on densities from the corresponding classes  $D_0^{k}$, $k=1,2,3,4$. The minimax spectral characteristic of the optimal estimate of the functional $A_s\vec{\xi}$ is determined by the formula (\ref{spchar2}).
\end{nas}

\begin{nas}
 The least favorable spectral densities $F^{0}(\lambda)$ in the classes $D_\varepsilon^{k}$, $k=1,2,3,4$, for the optimal linear interpolation of the functional  $A_s\vec{\xi}$, which depends on the unknown values of the process  $\vec{\xi}(t)$ based on observations of the process $\vec{\xi}(t)$ at points $t\in\mathbb{R} \backslash S$, are determined by the following  equations, respectively,
\begin{equation}
((C_{s}^{0}(\lambda) )^{\top} )^{*}\cdot(C_{s}^{0}(\lambda) )^{\top}=(\beta^{2}+\gamma(\lambda ))(F^{0} (\lambda ))^{2},
\end{equation}
\begin{equation}
((C_{s}^{0}(\lambda) )^{\top} )^{*}\cdot(C_{s}^{0}(\lambda) )^{\top}=F^{0} (\lambda )\left\{(\beta _{k}^{2} +\gamma_{k} (\lambda ))\delta _{kl} \right\}_{k,l=1}^{T}  F^{0} (\lambda ),
\end{equation}
\begin{equation}
((C_{s}^{0}(\lambda) )^{\top} )^{*}\cdot(C_{s}^{0}(\lambda) )^{\top}=(\beta^{2}+\gamma' (\lambda)) F^{0} (\lambda )(B_2)^\top F^{0} (\lambda ),
\end{equation}
\begin{equation}
((C_{s}^{0}(\lambda) )^{\top} )^{*}\cdot(C_{s}^{0}(\lambda) )^{\top}=F^{0} (\lambda ) (\vec{\beta}\cdot \vec{\beta}^{*}+\Gamma(\lambda ))F^{0} (\lambda ),
\end{equation}
 the minimality condition (\ref{minimal1}); the constrained optimization problem (\ref{extrem2}) and restrictions  on densities from the corresponding classes  $D_\varepsilon^{k}$, $k=1,2,3,4$. The minimax spectral characteristic of the optimal estimate of the functional $A_s\vec{\xi}$ is determined by the formula (\ref{spchar2}).
\end{nas}

\section{Least favorable spectral densities in the class $D=D_{V}^U \times D_{2\delta}$}

Consider the problem of the  mean-square optimal interpolation of the functional $A_s\vec{\xi}$  in the case where the spectral density matrices of the processes are not exactly known while the
 admissible spectral density matrices are from the class $D=D_{V}^U \times D_{2\delta}$, where
$${D_{V}^{U}} ^{1}  =\bigg\{F(\lambda )\bigg|{\mathrm{Tr}}\, V(\lambda
)\le {\mathrm{Tr}}\, F(\lambda )\le {\mathrm{Tr}}\, U(\lambda ), \frac{1}{2\pi } \int _{-\infty }^{\infty }{\mathrm{Tr}}\,  F(\lambda)d\lambda  =p \bigg\},$$
$$D_{2\delta}^{1}=\left\{G(\lambda )\biggl|\frac{1}{2\pi } \int_{-\infty }^{\infty }\left|{\rm{Tr}}(G(\lambda )-G_{1} (\lambda))\right|^{2} d\lambda \le \delta\right\};$$
$${D_{V}^{U}} ^{2}  =\bigg\{F(\lambda )\bigg|v_{kk} (\lambda )  \le
f_{kk} (\lambda )\le u_{kk} (\lambda ), \frac{1}{2\pi } \int _{-\infty }^{\infty }f_{kk} (\lambda
)d\lambda  =p_{k} , k=\overline{1,T}\bigg\},$$
$$D_{2\delta}^{2}=\left\{G(\lambda )\biggl|\frac{1}{2\pi } \int_{-\infty }^{\infty }\left|g_{kk} (\lambda )-g_{kk}^{1} (\lambda)\right|^{2} d\lambda  \le \delta_{k}, k=\overline{1,T}\right\};$$
$${D_{V}^{U}} ^{3}  =\bigg\{F(\lambda )\bigg|\left\langle B_{1}
,V(\lambda )\right\rangle \le \left\langle B_{1} ,F(\lambda
)\right\rangle \le \left\langle B_{1} ,U(\lambda )\right\rangle,\frac{1}{2\pi }
\int _{-\infty }^{\infty }\left\langle B_{1},F(\lambda)\right\rangle d\lambda  =p\bigg\},$$
$$D_{2\delta}^{3}=\left\{G(\lambda )\biggl|\frac{1}{2\pi } \int_{-\infty }^{\infty }\left|\left\langle B_{2} ,G(\lambda )-G_{1}(\lambda )\right\rangle \right|^{2} d\lambda  \le \delta\right\};$$
$${D_{V}^{U}} ^{4}=\left\{F(\lambda )\bigg|V(\lambda )\le F(\lambda
)\le U(\lambda ), \frac{1}{2\pi } \int _{-\infty }^{\infty}F(\lambda )d\lambda=P\right\}.$$
$$D_{2\delta}^{4}=\left\{G(\lambda )\biggl|\frac{1}{2\pi } \int_{-\infty }^{\infty }\left|g_{ij} (\lambda )-g_{ij}^{1} (\lambda)\right|^{2} d\lambda  \le \delta_{i}^{j}, i,j=\overline{1,T}\right\},$$
Here the spectral density matrices $V( \lambda ),U( \lambda ),G_{1} ( \lambda )$ are known and fixed,
$\delta, p, \delta_k, p_k, k=\overline{1,T}$, $\delta_{i}^{j}, i,j=\overline{1,T}$, are given numbers, $P, B_1, B_2$ are given positive-definite Hermitian matrices.

The class $ D_V^U$ describes the ``strip'' model of stochastic processes while the class $D_{2\delta}$ describes ``$\delta$-neighborhood''\, model  in the space $L_2$ of a given bounded spectral density  $G_1(\lambda)$.

From the condition $0\in \partial \Delta _{D} (F^{0} ,G^{0} )$ we find the following equations which determine the least favourable spectral densities for these given sets of admissible spectral densities.

For the first pair ${D_{V}^{U}} ^{1}\times D_{2\delta}^{1}$ we have equations
\[((A_{s} (\lambda ))^{\top} G^{0} (\lambda )+(C_{s}^{0} (\lambda
))^{\top} )^{*}((A_{s} (\lambda ))^{\top} G^{0} (\lambda )+(C_{s}^{0}
(\lambda ))^{\top} )=\]
\begin{equation} \label{eq_5_1}
=(\alpha ^{2} +\gamma _{1}
(\lambda )+\gamma _{2} (\lambda ))(F^{0} (\lambda )+G^{0} (\lambda))^{2},
\end{equation}
\[((A_{s} (\lambda ))^{\top} F^{0} (\lambda )-(C_{s}^{0} (\lambda
))^{\top} )^{*}((A_{s} (\lambda ))^{\top} F^{0} (\lambda )-(C_{s}^{0}
(\lambda ))^{\top})=\]
\begin{equation} \label{eq_5_2}
=\beta ^{2} {\mathrm{Tr}}\, (G^{0}(\lambda )-G_{1} (\lambda ))(F^{0} (\lambda )+G^{0} (\lambda ))^{2},
\end{equation}
\begin{equation} \label{eq_5_3}
\frac{1}{2\pi } \int _{-\infty }^{\infty }\left|{\mathrm{Tr}}\, (G(\lambda )-G_{1} (\lambda ))\right|^{2} d\lambda  =\delta ,
\end{equation}
where $\alpha^{2}, \beta^{2}$ are Lagrange multipliers, $\gamma _{1} (\lambda )\le 0$ and $\gamma _{1} (\lambda )=0$ if ${\mathrm{Tr}}\,
F^{0} (\lambda )> {\mathrm{Tr}}\,  V(\lambda ),$ $\gamma _{2} (\lambda )\ge 0$ and $\gamma _{2} (\lambda )=0$ if $ {\mathrm{Tr}}\,F^{0}(\lambda )< {\mathrm{Tr}}\,  U(\lambda).$

For the second pair ${D_{V}^{U}} ^{2}\times D_{2\delta}^{2}$ we have equations
\[((A_{s} (\lambda ))^{\top} G^{0} (\lambda )+(C_{s}^{0} (\lambda
))^{\top} )^{*}((A_{s} (\lambda ))^{\top} G^{0} (\lambda )+(C_{s}^{0}
(\lambda ))^{\top} )=\]
\begin{equation}\label{eq_5_4}
=(F^{0} (\lambda)+G^{0} (\lambda ))\left\{(\alpha_{k}^{2} +\gamma _{1k} (\lambda )+\gamma _{2k}(\lambda ))\delta _{kl} \right\}_{k,l=1}^{T} (F^{0} (\lambda )+G^{0}(\lambda )),
\end{equation}
\[((A_{s} (\lambda ))^{\top} F^{0} (\lambda )-(C_{s}^{0} (\lambda
))^{\top} )^{*}((A_{s} (\lambda ))^{\top} F^{0} (\lambda )-(C_{s}^{0}
(\lambda ))^{\top})=\]
\begin{equation} \label{eq_5_5}
=(F^{0} (\lambda )+G^{0}(\lambda )) \left\{\beta _{k}^{2} (g_{kk}^{0} (\lambda )-g_{kk}^{1}(\lambda ))\delta _{kl} \right\}_{k,l=1}^{T}
(F^{0} (\lambda )+G^{0}(\lambda )),
\end{equation}
\begin{equation} \label{eq_5_6}
\frac{1}{2\pi } \int _{-\infty }^{\infty }\left|g_{kk} (\lambda )-g_{kk}^{1} (\lambda )\right|^{2} d\lambda  =\delta _{k},\; k=\overline{1,T},
\end{equation}
where $\alpha _{k}^{2}, \beta_{k}^{2}$ are Lagrange multipliers, $\gamma _{1k} (\lambda )\le 0$ and $\gamma _{1k} (\lambda )=0$ if $f_{kk}^{0} (\lambda )>v_{kk} (\lambda ),$ $\gamma _{2k} (\lambda )\ge 0$ and $\gamma _{2k} (\lambda )=0$ if $f_{kk}^{0} (\lambda )<u_{kk} (\lambda).$

For the third pair ${D_{V}^{U}} ^{3}\times D_{2\delta}^{3}$ we have equations
\[((A_{s} (\lambda ))^{\top} G^{0} (\lambda )+(C_{s}^{0} (\lambda
))^{\top} )^{*}((A_{s} (\lambda ))^{\top} G^{0} (\lambda )+(C_{s}^{0}
(\lambda ))^{\top} )=\]
\begin{equation}\label{eq_5_7}
=(\alpha ^{2} +\gamma'_{1} (\lambda )+\gamma'_{2} (\lambda
))(F^{0} (\lambda )+G^{0} (\lambda ))(B_{1})^{\top}(F^{0} (\lambda)+G^{0} (\lambda )),
\end{equation}
\[((A_{s} (\lambda ))^{\top} F^{0} (\lambda )-(C_{s}^{0} (\lambda
))^{\top} )^{*}((A_{s} (\lambda ))^{\top} F^{0} (\lambda )-(C_{s}^{0}
(\lambda ))^{\top})=\]
\begin{equation} \label{eq_5_8}
=\beta ^{2} \left\langle B_{2},G^{0} (\lambda )-G_{1} (\lambda )\right\rangle(F^{0} (\lambda )+G^{0} (\lambda ))^2,
\end{equation}
\begin{equation} \label{eq_5_9}
\frac{1}{2\pi } \int _{-\infty }^{\infty }\left|\left\langle B_{2} ,G(\lambda )-G_{1} (\lambda )\right\rangle \right|^2d\lambda  =\delta,
\end{equation}
where $\alpha^{2}, \beta^{2}$ are Lagrange multipliers, $\gamma'_{1}( \lambda )\le 0$ and $\gamma'_{1} ( \lambda )=0$ if $\langle B_{1},F^{0} ( \lambda \rangle > \langle B_{1},V( \lambda ) \rangle,$ $\gamma'_{2}( \lambda )\ge 0$ and $\gamma'_{2} ( \lambda )=0$ if $\langle
B_{1} ,F^{0} ( \lambda \rangle < \langle B_{1} ,U( \lambda ) \rangle.$

For the fourth pair ${D_{V}^{U}} ^{4}\times D_{2\delta}^{4}$ we have equations
\[((A_{s} (\lambda ))^{\top} G^{0} (\lambda )+(C_{s}^{0} (\lambda
))^{\top} )^{*}((A_{s} (\lambda ))^{\top} G^{0} (\lambda )+(C_{s}^{0}
(\lambda ))^{\top} )=\]
\begin{equation}\label{eq_5_10}
=(F^{0} (\lambda )+G^{0} (\lambda ))(\vec{\alpha }\cdot \vec{\alpha
}^{*}+\Gamma _{1} (\lambda )+\Gamma _{2} (\lambda ))(F^{0} (\lambda)+G^{0} (\lambda ))
\end{equation}
\[((A_{s} (\lambda ))^{\top} F^{0} (\lambda )-(C_{s}^{0} (\lambda
))^{\top} )^{*}((A_{s} (\lambda ))^{\top} F^{0} (\lambda )-(C_{s}^{0}
(\lambda ))^{\top})=\]
\begin{equation} \label{eq_5_11}
=(F^{0} (\lambda )+G^{0} (\lambda ))\left\{\beta _{ij}^{} (g_{ij}^{0} (\lambda )-g_{ij}^{1}
(\lambda ))\right\}_{i,j=1}^{T} (F^{0} (\lambda )+G^{0} (\lambda )),
\end{equation}
\begin{equation} \label{eq_5_12}
\frac{1}{2\pi } \int _{-\infty }^{\infty }\left|g_{ij} (\lambda )-g_{ij}^{1} (\lambda )\right|^{2} d\lambda  =\delta_{i}^{j},\; i,j=\overline{1,T}.
\end{equation}
where  $\vec{\alpha}, \beta _{ij}$ are Lagrange multipliers, $\Gamma _{1} (\lambda )\le 0$ and $\Gamma _{1} (\lambda )=0$ if $F^{0}(\lambda )>V(\lambda ),$ $
\Gamma _{2} (\lambda )\ge 0$ and $\Gamma _{2} (\lambda )=0$ if $F^{0}(\lambda )<U(\lambda ).$

The following theorem and corollaries hold true.

\begin{thm}
 The least favorable spectral densities  $F^0(\lambda)$, $G^0(\lambda)$  in the classes $D=D_{V}^U \times D_{2\delta}$ for the optimal linear interpolation of the functional $A_s\vec{\xi}$ are determined by relations
(\ref{eq_5_1}) -- (\ref{eq_5_3}) for the first pair ${D_{V}^{U}} ^{1}\times D_{2\delta}^{1}$ of sets of admissible spectral densities;
(\ref{eq_5_4}) -- (\ref{eq_5_6}) for the second pair ${D_{V}^{U}} ^{2}\times D_{2\delta}^{2}$ of sets of admissible spectral densities;
(\ref{eq_5_7}) -- (\ref{eq_5_9}) for the third pair ${D_{V}^{U}} ^{3}\times D_{2\delta}^{3}$ of sets of admissible spectral densities;
(\ref{eq_5_10}) -- (\ref{eq_5_12}) for the fourth pair ${D_{V}^{U}} ^{4}\times D_{2\delta}^{4}$ of sets of admissible spectral densities;
the minimality condition (\ref{minimal}); the constrained optimization problem (\ref{extrem}) and restrictions  on densities from the corresponding classes $D=D_{V}^U \times D_{2\delta}$.  The minimax-robust spectral characteristic of the optimal estimate of the functional $A_s\vec{\xi}$ is determined by the formula (\ref{4}).
\end{thm}

\begin{nas}
 The least favorable spectral densities $F^{0}(\lambda)$ in the classes ${D_{V}^{U}} ^{k}$, $k=1,2,3,4$, for the optimal linear estimation of the functional  $A_s\vec{\xi}$,  which depends on the unknown values of the process  $\vec{\xi}(t)$ based on observations of the process $\vec{\xi}(t)$ at points $t\in\mathbb{R} \backslash S$, are determined by the following  equations, respectively,
\begin{equation}
((C_{s}^{0}(\lambda) )^{\top} )^{*}\cdot(C_{s}^{0}(\lambda) )^{\top}=(\alpha^{2} +\gamma _{1} (\lambda )+\gamma _{2} (\lambda )) (F^{0} (\lambda ))^{2},
\end{equation}
\begin{equation}
((C_{s}^{0}(\lambda) )^{\top} )^{*}\cdot(C_{s}^{0}(\lambda) )^{\top}=F^{0} (\lambda )\left\{\beta _{k}^{2} (g_{kk}^{0} (\lambda )-g_{kk}^{1}(\lambda ))\delta _{kl} \right\}_{k,l=1}^{T}
F^{0} (\lambda ),
\end{equation}
\begin{equation}
((C_{s}^{0}(\lambda) )^{\top} )^{*}\cdot(C_{s}^{0}(\lambda) )^{\top}=(\alpha^{2} +\gamma'_{1}(\lambda )+\gamma'_{2}(\lambda )) F^{0} (\lambda )(B_1)^\top F^{0} (\lambda ),
\end{equation}
\begin{equation}
((C_{s}^{0}(\lambda) )^{\top} )^{*}\cdot(C_{s}^{0}(\lambda) )^{\top}=F^{0} (\lambda )(\vec{\alpha}\cdot \vec{\alpha}^{*}+\Gamma _{1} (\lambda )+\Gamma _{2} (\lambda ))F^{0} (\lambda ),
\end{equation}
the minimality condition (\ref{minimal}); the constrained optimization problem (\ref{extrem2}) and restrictions  on densities from the corresponding classes  ${D_{V}^{U}} ^{k}$, $k=1,2,3,4$. The minimax spectral characteristic of the optimal estimate of the functional $A_s\vec{\xi}$ is determined by the formula (\ref{spchar2}).
\end{nas}

\begin{nas}
 The least favorable spectral densities $F^{0}(\lambda)$ in the classes $D_{2\delta}^{k}$, $k=1,2,3,4$, for the optimal linear estimation of the functional  $A_s\vec{\xi}$,  which depends on the unknown values of the process  $\vec{\xi}(t)$ based on observations of the process $\vec{\xi}(t)$ at points $t\in\mathbb{R} \backslash S$, are determined by the following  equations, respectively,
\begin{equation}
((C_{s}^{0}(\lambda) )^{\top} )^{*}\cdot(C_{s}^{0}(\lambda) )^{\top}=\beta ^{2} {\mathrm{Tr}}\, (F^{0} (\lambda )-G_{1} (\lambda ))(F^{0} (\lambda))^{2} ,
\end{equation}
\begin{equation}
((C_{s}^{0}(\lambda) )^{\top} )^{*}\cdot(C_{s}^{0}(\lambda) )^{\top}=F^{0} (\lambda )\left\{\beta _{k}^{2} (f_{kk}^{0} (\lambda)-g_{kk}^{1} (\lambda ))\delta _{kl} \right\}_{k,l=1}^{T} F^{0}(\lambda ),
\end{equation}
\begin{equation}
((C_{s}^{0}(\lambda) )^{\top} )^{*}\cdot(C_{s}^{0}(\lambda) )^{\top}=\beta ^{2} \left\langle B_{2} ,F^{0} (\lambda )-G_{1} (\lambda)\right\rangle (F^{0} (\lambda ))^2,
\end{equation}
\begin{equation}
((C_{s}^{0}(\lambda) )^{\top} )^{*}\cdot(C_{s}^{0}(\lambda) )^{\top}=F^{0} (\lambda )\left\{\beta _{ij}^{} (f_{ij}^{0} (\lambda)-g_{ij}^{1} (\lambda ))\right\}_{i,j=1}^{T} F^{0} (\lambda )
\end{equation}
the minimality condition (\ref{minimal}); the constrained optimization problem (\ref{extrem2}) and the following restrictions  on densities from the corresponding classes  $D_{2\delta}^{k}$, $k=1,2,3,4$, respectively,
\begin{equation}
\frac{1}{2\pi } \int _{-\infty }^{\infty }\left|{\mathrm{Tr}}\, (G(\lambda )-G_{1} (\lambda ))\right|^{2} d\lambda  =\delta ,
\end{equation}
\begin{equation}
\frac{1}{2\pi } \int _{-\infty }^{\infty }\left|g_{kk} (\lambda )-g_{kk}^{1} (\lambda )\right|^{2} d\lambda  =\delta _{k},\; k=\overline{1,T},
\end{equation}
\begin{equation}
\frac{1}{2\pi } \int _{-\infty }^{\infty }\left|\left\langle B_{2} ,G(\lambda )-G_{1} (\lambda )\right\rangle \right|^2d\lambda  =\delta,
\end{equation}
\begin{equation}
\frac{1}{2\pi } \int _{-\infty }^{\infty }\left|g_{ij} (\lambda )-g_{ij}^{1} (\lambda )\right|^{2} d\lambda  =\delta_{i}^{j},\; i,j=\overline{1,T}.
\end{equation}
The minimax spectral characteristic of the optimal estimate of the functional $A_s\vec{\xi}$ is determined by the formula (\ref{spchar2}).
\end{nas}

\section{Conclusions}
In this article we propose methods of the mean-square optimal linear interpolation of functionals which depend on the unknown values of a multidimensional continuous time  stationary stochastic process based on observations of the process with an additive stationary stochastic noise process.
 The case of spectral certainty as well as the case  of spectral uncertainty are considered.
Under the condition of spectral certainty, where the spectral density matrices of the stationary processes are exactly known, we derive formulas for calculating the spectral characteristic and the mean-square error of the optimal estimate of the functional.
The corresponding results are derived in the case of observations of the process without noise.
In the case of spectral uncertainty, where the spectral density matrices of the stationary processes are not exactly known
while some sets of admissible spectral density matrices are given, we apply the minimax-robust method of estimation.
This method allows us to find estimates that minimize the maximum values of the mean-square errors of the estimates for all spectral density matrices from a given class of admissible spectral density matrices
and derive relations which determine the least favourable spectral density matrices.
These least favourable spectral density matrices are solutions of the optimization problem (\ref{8})
which is characterized by the condition $0 \in \partial\Delta_D(F^0,G^0),$ where $\partial\Delta_D(F^0,G^0)$ is the subdifferential of the convex functional $\Delta_D(F,G)$ at point $(F^0,G^0)$.
The form of the functional $\Delta(h(F^0,G^0);F,G)$ is convenient for application of the Lagrange method of indefinite multipliers for finding solution to the problem (\ref{8}).
Making use of the method of Lagrange multipliers and the form of
subdifferentials of the indicator functions we describe relations that determine the least favourable spectral densities in some special classes of spectral densities.

\end{document}